# Phase-only signal reconstruction by MagnitudeCut

Jiasong Wu[1,2,3,4], Jieyuan Liu[1,4], Youyong Kong[1,4], Xu Han[1,4], Lotfi Senhadji[2,3,4], Huazhong Shu[1,4]

[1]*LIST, Key Laboratory of Computer Network and Information Integration (Southeast University), Ministry of Education, Nanjing, China*

[2]*INSERM, U1099, Rennes, France*

[3]*LTSI, Université de Rennes 1, Rennes, France*

[4] *Centre de Recherche en Information Biomédicale Sino-Français (CRIBs)*

E-mail*:* jswu@seu.edu.cn, 1095486276@qq.com, kongyouyong@seu.edu.cn, 229897099@qq.com, lotfi.senhadji@univ-rennes1.fr, shu.list@seu.edu.cn

Corresponding author: Jiasong Wu ( *jswu@seu.edu.cn*)

The authors declare that there is no conflict of interest regarding the publication of this paper.

**Abstract:** In this paper, we present a new algorithm, called MagnitudeCut, for recovering a signal from the phase of its Fourier transform. We casted our recovering problem into a new convex optimization problem, and then solved it by the block coordinate descent algorithm and the interior point algorithm, in which the iteration process consists of matrix vector product and inner product. We used the new method for reconstruction of a set of signal/image. The simulation results reveal that the proposed MagnitudeCut method can reconstruct the original signal with fewer sampling number of the phase information than that of the Greedy algorithm and iterative method under the same reconstruction error. Moreover, our algorithm can also reconstruct the symmetric image from its Fourier phase.

## 1. Introduction

Generally speaking, the Fourier phase and the Fourier magnitude of a signal are mutually independent, so the signal cannot be recovered only from the knowledge of one of them. However, Hayes and Oppenheim [1] first pointed out that it is possible to recover the original signal from the phase-only information under certain conditions: when a finite length sequence has a *z*-transform with no zero on the unit circle or in conjugate reciprocal pairs, it can be uniquely specified by a scale factor based on some distinct samples of its Fourier phase. Subsequently, Levi and Stark [2] suggested a new condition, that is, the *z*-transform of the signal cannot contain symmetric factors, in other words, the original signal has no symmetric or antisymmetric axis. Recently, Boufounos et al. [3] exploited sparse signal models in order to reduce the number of Fourier phase measurements required in the reconstruction.

Since then, the problem of restoring a signal from the Fourier phase continues to receive great attention. During the past decades, many algorithms have been reported in the literature to solve the above problem. Among them, there are iterative methods and closed-form solutions [1, 4], the alternating projection methods [2], and partial phase information methods [5, 6]. Hua and Orchard [7-9] proposed a new image reconstruction algorithm with simple geometrical models. Loveimi and Ahadi [10] reconstructed the speech signal via the least square error estimation and the overlap-add methods. Recently, Boufounos [3] realized reconstruction by leveraging standard convex and greedy algorithms. It should be pointed out that in all the aforementioned methods, the conditions were imposed on the original signal, and they are not suitable to recover the symmetrical signal from the Fourier phase.

In this paper, inspired by the recent work of PhaseLift [11] and PhaseCut [12] for solving the phase retrieval problem, we deal with the corresponding magnitude retrieval problem, that is, signal/image reconstruction from its Fourier phase information. We follow the way reported in [11, 12]. That is, we do not restrict the conditions on the original signal, but consider the sampling number of its Fourier phase. To capture sufficient Fourier phase information, we used one or multiple random masks to multiply by the original signal before Fourier transform. In our proposed method (called MagnitudeCut method), we used the block coordinate descent algorithm and the interior point method during the iteration process, due to their simplicity and efficiency.

We solve the signal reconstruction problem by exploiting only the phase information because it is well known that most of the signal information contained in the phase is more important than that contained in the magnitude under the same number of samples [13]. Moreover, the phase has been found very important in many applications such as image retrieval [14], object recognition [15]. Hence, it can be expected that by taking the advantage of phase information, the signal can be reconstructed with less samples than that based on magnitude.

The paper is organized as follows. The MagnitudeCut method is proposed in section 2. Section 3 gives the experiment results on simulation data and also real data. Section 4 concludes the work and gives the future work.

## 2. Phase-only signal reconstruction

In this section, we first derived the MagnitudeCut method for recovering the complex signal $\boldsymbol{x} \in \mathbb{C}^p$ by using the block coordinate descent algorithm and log-barrier algorithm in subsection 2.1, then we discussed why our proposed algorithm can deal with the symmetrical signal which can not be solved by the algorithms in [1, 2, 13] in subsection 2.2. In subsection 2.3 we give some remarks on the proposed MagnitudeCut method.

The original problem is formulated as follow:

$$\begin{aligned} &\text{find} \quad \boldsymbol{x} \\ &\text{such that} \quad e^{j\angle \boldsymbol{Ax}} = \boldsymbol{u}, \end{aligned} \tag{1}$$

where $\boldsymbol{x} \in \mathbb{C}^p$ is the signal that we want to recovery, and $\mathbb{C}$ denotes the complex number domain; $\angle \boldsymbol{Ax}$ denotes the phase of $\boldsymbol{Ax}$, where the transform matrix $\boldsymbol{A} \in \mathbb{C}^{n\times p}$ is the product of Fourier matrix $\boldsymbol{F} \in \mathbb{C}^{p\times p}$ and random mask matrix $\boldsymbol{R}_i \in \mathbb{C}^{p\times p}$, $i = (1,2,...,n/p)$, i.e. $\boldsymbol{A}=\begin{bmatrix} \boldsymbol{FR}_1 & \boldsymbol{FR}_2 & \cdots & \boldsymbol{FR}_{n/p} \end{bmatrix}^T$ and $n$ is a multiple of $p$; $\boldsymbol{u} \in \mathbb{C}^n$ is the value that we obtained and its elements $u_i, i = 1,...,n$, satisfy $|u_i| = 1$.

The advantage of the transform matrix $\boldsymbol{A}$ is two-fold. First, enough Fourier phase information can be obtained by multiplying the signal with multiple random masks. Second, the original symmetric signal usually becomes asymmetric after multiplying by a random mask. Note that the idea of multiplying the random masks comes from [11, 12], readers may refer to the mentioned two references for more details.

### *2.1 MagnitudeCut Method*

As mentioned above, the transform matrix $\boldsymbol{A}$ is composed by random mask matrix and Fourier matrix. We used different random masks for obtaining sufficient phase information. Our objective is to recover the original signal with an appropriate number of phase information. We separate the Fourier magnitude from Fourier phase variables in order that the representation of signal in Fourier domain can be written as $\boldsymbol{Ax} = \text{diag}(\boldsymbol{u})\boldsymbol{b}$, where $\boldsymbol{b} \in \mathbb{R}^n$ denotes the Fourier magnitude and $\mathbb{R}$ denotes the real number domain and $\text{diag}(\boldsymbol{u}) = \begin{bmatrix} u_1 & & & \\ & u_2 & & \\ & & \ddots & \\ & & & u_n \end{bmatrix}$. Given $e^{j\angle \boldsymbol{Ax}} = \boldsymbol{u}$, (1) can thus be written as a least square problem

$$\min_{\boldsymbol{x}\in\mathbb{C}^p,\boldsymbol{b}\in\mathbb{R}^n} \left\| \boldsymbol{Ax} - \text{diag}(\boldsymbol{u})\boldsymbol{b} \right\|_2^2, \tag{2}$$

where $\|\bullet\|_2$ is the $\ell_2$-norm and defined as $\|\bullet\|_2 = \left( \sum_{i=1}^{n} |\bullet|^2 \right)^{\frac{1}{2}}$, that is, $\|\bullet\|_2^2 = \sum_{i=1}^{n} |\bullet|^2$. By substituting $\boldsymbol{x}$ in (2) with $\boldsymbol{x} = \boldsymbol{A}^\dagger \text{diag}(\boldsymbol{u})\boldsymbol{b}$, where $(\bullet)^\dagger$ is the pseudo-inverse operator, then (2) with respect to $\boldsymbol{x}$ and $\boldsymbol{b}$ can be converted to a problem only relative to $\boldsymbol{b}$. It means that (2) is equivalent to

$$\min_{\boldsymbol{b}\in\mathbb{R}^n} \left\| \boldsymbol{AA}^\dagger \text{diag}(\boldsymbol{u})\boldsymbol{b} - \text{diag}(\boldsymbol{u})\boldsymbol{b} \right\|_2^2. \tag{3}$$

The objective function of this problem can be rewritten as

$$\left\| \boldsymbol{AA}^\dagger \text{diag}(\boldsymbol{u})\boldsymbol{b} - \text{diag}(\boldsymbol{u})\boldsymbol{b} \right\|_2^2 = \boldsymbol{b}^T \text{diag}(\boldsymbol{u}^H) \tilde{\boldsymbol{M}} \text{diag}(\boldsymbol{u})\boldsymbol{b},$$

where $\tilde{\boldsymbol{M}} = (\boldsymbol{AA}^\dagger - \boldsymbol{I})^H (\boldsymbol{AA}^\dagger - \boldsymbol{I}) = \boldsymbol{I} - \boldsymbol{AA}^\dagger$, the superscript $H$ is the conjugate transpose operator and $\boldsymbol{I}$ is the identity matrix. Let $\boldsymbol{M}$ be a positive semi-definite matrix given by $\boldsymbol{M} = \text{diag}(\boldsymbol{u}^H)(\boldsymbol{I} - \boldsymbol{AA}^\dagger)\text{diag}(\boldsymbol{u})$, so that (3) becomes

$$\min_{\boldsymbol{b}} \boldsymbol{b}^T \boldsymbol{Mb}, \quad s.t. \quad \boldsymbol{b} \in \mathbb{R}^n. \tag{4}$$

Next, letting $\boldsymbol{B} = \boldsymbol{bb}^T \in \mathrm{S}_n^+$, where $\mathrm{S}_n^+$ denotes the set of positive symmetric matrices, then (4) becomes

$$\begin{aligned} &\min_{\boldsymbol{B}} \quad \text{Tr}(\boldsymbol{BM}) \\ &s.t. \quad \boldsymbol{B} \geq 0, \boldsymbol{B} \in S_n^+, \text{rank}(\boldsymbol{B}) = 1. \end{aligned} \tag{5}$$

After dropping the non-convex rank constraint [12], we obtain the following convex relaxation problem

$$\min_{\boldsymbol{B}} \text{Tr}(\boldsymbol{BM}), \quad s.t. \quad \boldsymbol{B} \geq 0. \tag{6}$$

Since both $\boldsymbol{M}$ and $\boldsymbol{B}$ belong to $\mathrm{H}_n$, where $\mathrm{H}_n$ is the Hermitian matrices of dimension $n$, we can formulate the complex operation in $\mathrm{H}_n$ as the real operation [16]. Given the function $\Gamma(\cdot) = \begin{bmatrix} \text{Re}(\cdot) & -\text{Im}(\cdot) \\ \text{Im}(\cdot) & \text{Re}(\cdot) \end{bmatrix}$ as in [17], we have

$$\text{Tr}(\boldsymbol{BM}) = \frac{1}{2} \text{Tr}(\Gamma(\boldsymbol{B})\Gamma(\boldsymbol{M})). \tag{7}$$

Applying $\Gamma(\cdot)$ to (7), we have

$$\text{Tr}(\Gamma(\boldsymbol{B})\Gamma(\boldsymbol{M})) = \text{Tr}(2(\boldsymbol{B}\text{Re}(\boldsymbol{M}))). \tag{8}$$

Finally, (6) becomes

$$\min_{\boldsymbol{B}} \text{Tr}(\boldsymbol{B}\text{Re}(\boldsymbol{M})), \quad \text{s.t.} \quad \boldsymbol{B} \geq 0. \tag{9}$$

*2.1.1. Block coordinate descent algorithm*

The proposed MagnitudeCut method is applied to the barrier version of the MaxCut [18] to relax the matrix $\boldsymbol{B}$, then the convex optimization problem (9) becomes

$$\min_{\boldsymbol{B}} \text{Tr}(\boldsymbol{B}\text{Re}(\boldsymbol{M})) - \mu \log(\det(\boldsymbol{B})), \quad \mu > 0, \tag{10}$$

where $-\log(\det(\boldsymbol{B}))$ is the log-barrier term and $\det(\boldsymbol{B})$ is the determinant of $\boldsymbol{B}$. In order to work out the minimum value of the objective function in (10), firstly we use the block coordinate descent (BCD) algorithm [16] to simplify the issue (10) so as to solve it more conveniently. We will introduce the reduction process below. We first rewrite the matrix $\boldsymbol{B}=\begin{bmatrix}\boldsymbol{P} & \boldsymbol{y}\\ \boldsymbol{y}^T & b_n^2\end{bmatrix}$, where $\boldsymbol{P}=\begin{bmatrix}b_1^2 & \cdots & b_1b_{n-1}\\ \vdots & \ddots & \vdots\\ b_{n-1}b_1 & \cdots & b_{n-1}^2\end{bmatrix}$,. $\boldsymbol{y}=\begin{bmatrix}b_1b_n\\ \vdots\\ b_{n-1}b_n\end{bmatrix}=b_n\begin{bmatrix}b_1\\ \vdots\\ b_{n-1}\end{bmatrix}=b_n\boldsymbol{y}'$; then according to [19], we know

$$\begin{bmatrix}\boldsymbol{P} & \boldsymbol{y}\\ \boldsymbol{y}^T & b^2\end{bmatrix}\begin{bmatrix}\boldsymbol{I} & -\boldsymbol{P}^{-1}\boldsymbol{y}\\ 0 & \boldsymbol{I}\end{bmatrix}=\begin{bmatrix}\boldsymbol{P} & 0\\ \boldsymbol{y}^T & b^2-\boldsymbol{y}^T\boldsymbol{P}^{-1}\boldsymbol{y}\end{bmatrix}.$$

As a consequence, $\det(\boldsymbol{B})=\det(\boldsymbol{P})\det(b^2-\boldsymbol{y}^T\boldsymbol{P}^{-1}\boldsymbol{y})$. Without loss of generality, we redescribe the determinant by using the elementary matrix transformation:

$$\det(\boldsymbol{B})=\det(\boldsymbol{P}_{i^ci^c})\det(b_i^2-\boldsymbol{y}_{i^ci}^T\boldsymbol{P}_{i^ci^c}^{-1}\boldsymbol{y}_{i^ci}), i=\{1,\ldots,n\}, i^c\in\{1,\ldots,i-1,i+1,\ldots,n\}. \tag{11}$$

Based on the derivation of above formulations, the minimization of (10) is equivalent to

$$\min_{\boldsymbol{y}'_{i^ci},b_i}\sum_i\{\mathrm{Re}(\boldsymbol{M})_{i^ci}^T\boldsymbol{y}'_{i^ci}b_i+\mathrm{Re}(\boldsymbol{M})_{ii}b_i^2\}-\mu\log(1-\boldsymbol{y}'^T_{i^ci}\boldsymbol{P}_{i^ci^c}^{-1}\boldsymbol{y}'_{i^ci})-\mu\log b_i^2. \tag{12}$$

Ultimately, (12) can be cast as the following block one,

$$\min_{\boldsymbol{y}'_{i^ci},b_i}\mathrm{Re}(\boldsymbol{M})_{i^ci}^T\boldsymbol{y}'_{i^ci}b_i+\mathrm{Re}(\boldsymbol{M})_{ii}b_i^2-\mu\log(1-\boldsymbol{y}'^T_{i^ci}\boldsymbol{P}_{i^ci^c}^{-1}\boldsymbol{y}'_{i^ci})-\mu\log b_i^2,\ i=1,\ldots,n. \tag{13}$$

*2.1.2. Log-barrier algorithm*

In the following context, we use the Log-barrier (LB) algorithm, which is the one of the classical interior point algorithms [20], to solve the problem (13). For convenience, we set $\boldsymbol{y}'=\boldsymbol{y}'_{i^ci}, \boldsymbol{Q}=\mathrm{Re}(\boldsymbol{M})_{i^ci}, q=\mathrm{Re}(\boldsymbol{M})_{ii}$ and $\boldsymbol{P}=\boldsymbol{P}_{i^ci^c}$, then the objective function in (13) can be reduced to $F(\boldsymbol{y}',b_i)=\boldsymbol{Q}^T\boldsymbol{y}'b_i+qb_i^2-\mu\log(1-\boldsymbol{y}'^T\boldsymbol{P}^{-1}\boldsymbol{y}')-\mu\log b_i^2$. Taking the derivative of $F(\boldsymbol{y}',b_i)$ with respect to $\boldsymbol{y}'$ and $b_i$, we get

$$\partial F(\boldsymbol{y}',b_i)/\partial\boldsymbol{y}'=\boldsymbol{Q}b_i+2\mu\boldsymbol{P}^{-1}\boldsymbol{y}'/(1-\boldsymbol{y}'^T\boldsymbol{P}^{-1}\boldsymbol{y}')=I(\boldsymbol{y}'), \tag{14}$$

$$\partial F(\boldsymbol{y}',b_i)/\partial b_i=\boldsymbol{Q}^T\boldsymbol{y}'+2qb_i-2\mu/b_i=J(b_i). \tag{15}$$

Next, by using the second order Taylor expansion, (14) and (15) can be expressed as

$$I(\boldsymbol{y}'+\Delta\boldsymbol{y}')=I(\boldsymbol{y}')+(2\mu(\boldsymbol{P}^{-1})^T/(1-\boldsymbol{y}'^T\boldsymbol{P}^{-1}\boldsymbol{y}'))\Delta\boldsymbol{y}'+(4\mu\boldsymbol{P}^{-1}\boldsymbol{y}'(\boldsymbol{P}^{-1}\boldsymbol{y}')^T/(1-\boldsymbol{y}'^T\boldsymbol{P}^{-1}\boldsymbol{y}')^2)\Delta\boldsymbol{y}', \tag{16}$$

$$J(b_i+\Delta b_i)=J(b_i)+2q\Delta b_i+(2\mu/b_i^2)\Delta b_i. \tag{17}$$

To simplify (16) and (17), we define

$$\begin{aligned}\boldsymbol{G}&=-\mu\log(1-\boldsymbol{y}'^T\boldsymbol{P}^{-1}\boldsymbol{y}')\\ \boldsymbol{g}_{y'}&=d\boldsymbol{G}/d\boldsymbol{y}'=2\mu\boldsymbol{P}^{-1}\boldsymbol{y}'/(1-\boldsymbol{y}'^T\boldsymbol{P}^{-1}\boldsymbol{y}')\\ \boldsymbol{H}_{y'}&=d^2\boldsymbol{G}/d\boldsymbol{y}'^2=2\mu(\boldsymbol{P}^{-1})^T/(1-\boldsymbol{y}'^T\boldsymbol{P}^{-1}\boldsymbol{y}')+4\mu\boldsymbol{P}^{-1}\boldsymbol{y}'(\boldsymbol{P}^{-1}\boldsymbol{y}')^T/(1-\boldsymbol{y}'^T\boldsymbol{P}^{-1}\boldsymbol{y}')^2\end{aligned}, \tag{18}$$

$$K=qb_i^2-\mu\log(b_i^2),\ k_{b_i}=2qb_i-2\mu/b_i,\ L_{b_i}=2q+2\mu/b_i^2. \tag{19}$$

Hence, (16) and (17) are simplified as follows

$$\begin{aligned}\boldsymbol{Q}b_i+\boldsymbol{g}_{y'}+\boldsymbol{H}_{y'}\Delta\boldsymbol{y}'&=\boldsymbol{0}\\ \boldsymbol{Q}^T\boldsymbol{y}'+k_{b_i}+L_{b_i}\Delta b_i&=0\end{aligned}. \tag{20}$$

The solution of (20) can be written in the following form

$$\begin{aligned}\Delta\boldsymbol{y}'&=-(\boldsymbol{H}_{y'})^{-1}(\boldsymbol{Q}b_i+\boldsymbol{g}_{y'})\\ \Delta b_i&=-(\boldsymbol{Q}^T\boldsymbol{y}'+k_{b_i})/L_{b_i}\end{aligned}. \tag{21}$$

The iterative equations are given as follows

$$\begin{aligned}\boldsymbol{y}'_{i^ci}&=\boldsymbol{y}'_{i^ci}+s\Delta\boldsymbol{y}'_{i^ci}\\ b_i&=b_i+s\Delta b_i\end{aligned}. \tag{22}$$

where *s* is the step size.

We update $\boldsymbol{y}'$ and $b_i$ by (21), so as to get a better solution $\hat{\boldsymbol{B}}$. If $\hat{\boldsymbol{B}}$ has rank one, the relaxation is tight and the vector $\hat{\boldsymbol{b}}$ such that $\hat{\boldsymbol{B}}=\hat{\boldsymbol{b}}\hat{\boldsymbol{b}}^T$ is an optimal solution of (4). When the rank of the solution $\hat{\boldsymbol{B}}$ is larger than one, $\sqrt{\lambda_{\max}}\nu_{\max}$ is used as an approximate solution [12], where $\lambda_{\max}$ is the maximum eigenvalue and $\nu_{\max}$ is the corresponding eigenvector. At last, we obtain the reconstruction signal by $\hat{\boldsymbol{x}} = \boldsymbol{A}^{\dagger}\mathrm{diag}(\boldsymbol{u})\hat{\boldsymbol{b}}$.

*2.2. A special situation*

As we have mentioned in the introduction, the traditional methods are not suitable to reconstruct a symmetric signal. In what follows, we will prove that the reconstruction of a symmetrical signal can be realized by MagnitudeCut method.

For a symmetric signal $\boldsymbol{x}_{(sym)}$ , according to the above context, we know

$$\boldsymbol{x}_{(asym)} = \boldsymbol{R}\boldsymbol{x}_{(sym)}, \tag{23}$$

$$\boldsymbol{X} = \boldsymbol{F}\boldsymbol{x}_{(asym)} = \boldsymbol{F}\boldsymbol{R}\boldsymbol{x}_{(sym)} = (\boldsymbol{F}\boldsymbol{R})\boldsymbol{x}_{(sym)} = \boldsymbol{A}\boldsymbol{x}_{(sym)}, \tag{24}$$

where $\boldsymbol{x}_{(asym)}$ is the asymmetric signal. The random mask $\boldsymbol{R}$ converts the symmetrical signal to an asymmetric signal. When we obtain the presentation of the Fourier domain $\boldsymbol{X}$ by MagnitudeCut, we can reconstruct the original symmetric signal through the following process.

$$\boldsymbol{x}_{(asym)} = \boldsymbol{F}^{-1}\boldsymbol{X}, \tag{25}$$

$$\boldsymbol{x}_{(sym)} = (\boldsymbol{A})^{-1}\boldsymbol{X} = (\boldsymbol{F}\boldsymbol{R})^{-1}\boldsymbol{X} = \boldsymbol{R}^{-1}\boldsymbol{F}^{-1}\boldsymbol{X} = \boldsymbol{R}^{-1}(\boldsymbol{F}^{-1}\boldsymbol{X}) = \boldsymbol{R}^{-1}\boldsymbol{x}_{(asym)}. \tag{26}$$

Up to now, we prove that MagnitudeCut can reconstruct the symmetric signal in theory. We provide some simulation experiments of symmetric signal reconstruction in section 3.

*2.3. Discussion*

In this subsection, we give some remarks about the proposed MagnitudeCut method. After obtaining $\hat{\boldsymbol{x}}$, we can set $\boldsymbol{X}^0 = \boldsymbol{A}\hat{\boldsymbol{x}} \in \mathrm{F}$ (the set F which satisfies (1)) as the initial value of the iteration method [1, 4], to get a more accurate solution. We know that the MagnitudeCut method needs only real matrix vector product and real inner product in each iteration process. Although we get the result by iteration method at last, the iteration method has no significant contribution to the computational complexity. The reason is that $\hat{\boldsymbol{x}}$ is sufficiently close to $\boldsymbol{x}$, thus the iteration method requires much less arithmetic operations than MagnitudeCut method. Since the critical component of the MagnitudeCut method is interior point algorithm, so the convergence of the algorithm is guaranteed by the result in [20].

## 3. Experiments on simulation data and real data

In this section, we first use the MagnitudeCut method to the reconstruction of one-dimensional (1-D) complex simulation signal in subsection 3.1. The reconstruction of two-dimensional (2-D) real medical images and 2-D complex symmetric images are given in subsections 3.2 and 3.3, respectively. In subsection 3.4, we compare the MagnitudeCut method with Greedy algorithm [23] and Iterative algorithm [19].

The following experiments were implemented in Matlab programming language on a PC machine, which sets up Microsoft Windows 7 operating system and has an Intel(R) Core(TM) i5-2400 CPU with speed of 3.10GHz and 4GB RAM.

*3.1. 1-D signal reconstruction*

The original signals $\boldsymbol{x} \in \mathbb{C}^{32}$ are shown in Fig. 1(a)-(d), and suppose the phase of its Fourier transform is also given. The sampling number is equal to the length of original signal, and the reconstruction results are shown in Fig. 1(A)-(D). From this figure, we find that the reconstructed results by MagnitudeCut method are perfect. In addition, from Fig. 1(c), (d), (C) and (D), we can see that the proposed MagnitudeCut method can also reconstruct the symmetric signal from its Fourier phase.

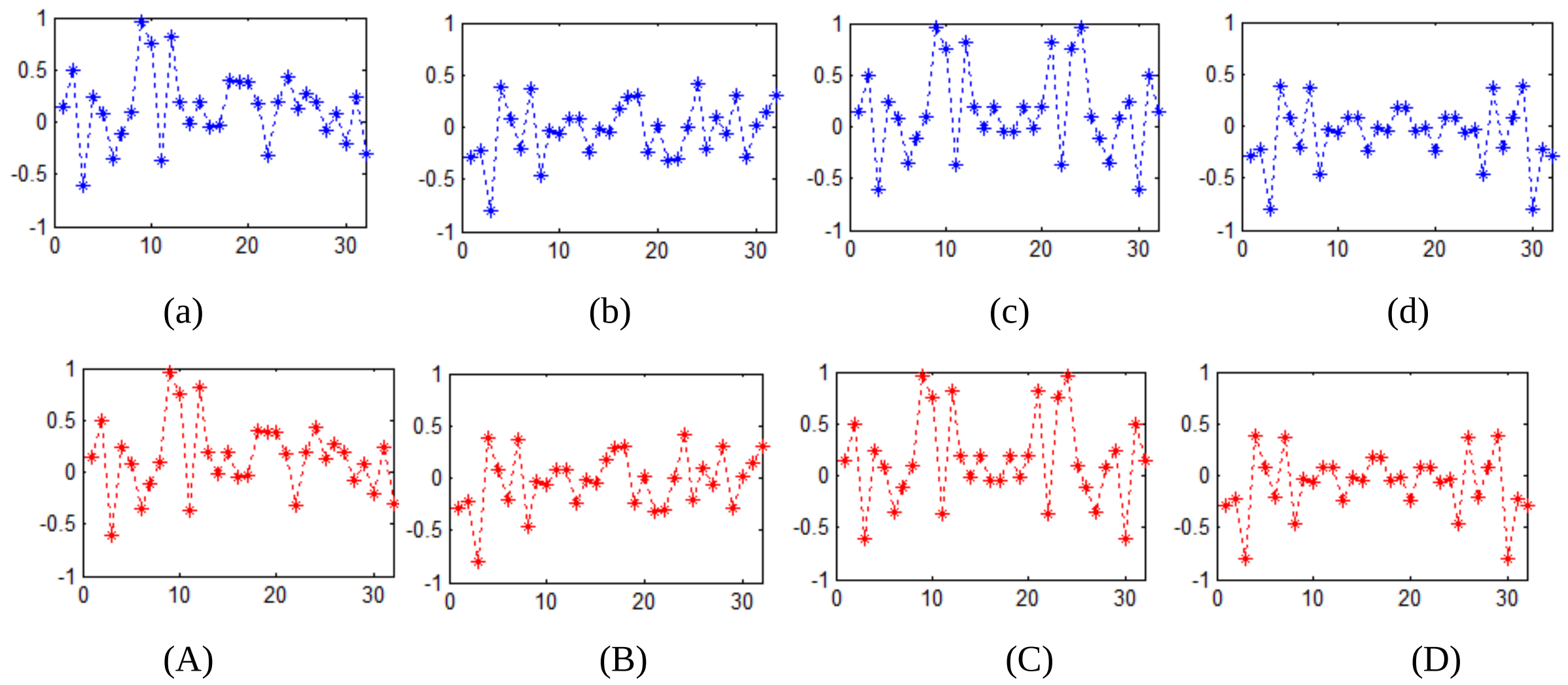

Fig. 1 The reconstruction results of the 1D complex asymmetric (left) and symmetric (right) signal. (a) and (c) are the real part of the original signals, (b) and (d) are the imaginary part of the original signals. (A) and (C) are the real part of the reconstruction signals, (B) and (D) are the imaginary part of the reconstruction signals.

*3.2. 2-D medical images reconstruction*

The medical images $\boldsymbol{x} \in \mathbb{R} \times \mathbb{R}$ are shown in the top line of Fig. 2 (a)-(f) from left to light, respectively. We suppose that the phases of their Fourier transform are also given. The *s_n* , which is the number of samples of phase and we call it sampling number in the following, is triple the length of original signal and the reconstruction results $\hat{\boldsymbol{x}}$ are shown in the bottom of Fig. 2.

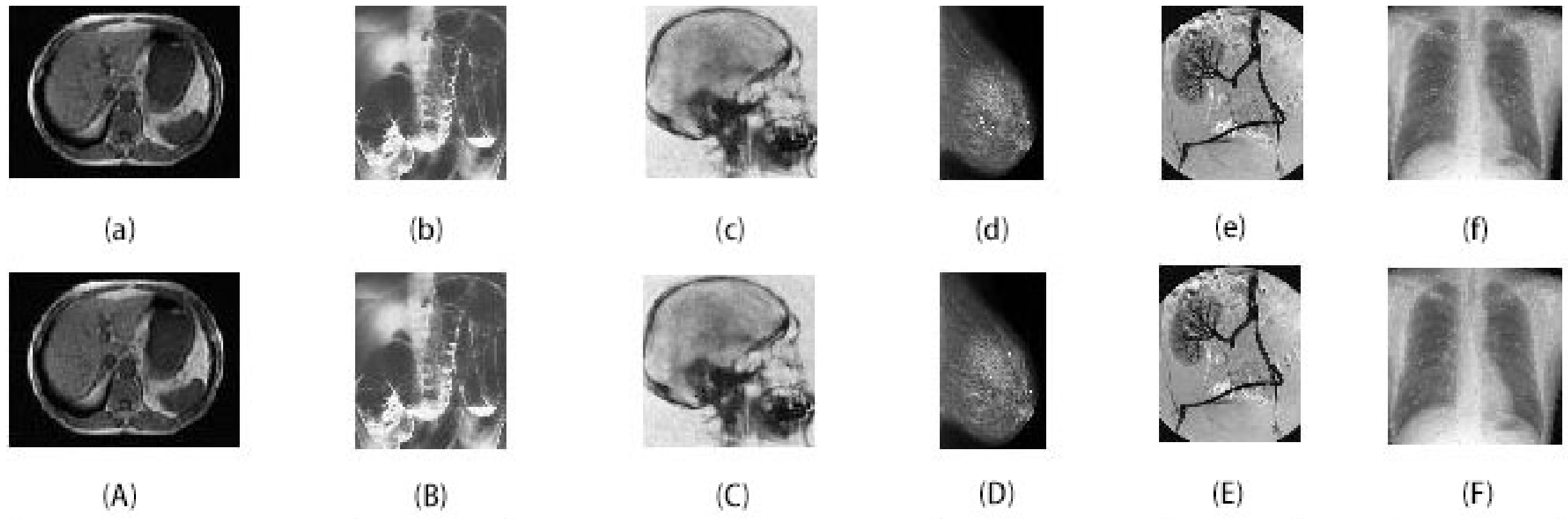

Fig. 2 (a)-(f) are original images, (A)-(F) are the reconstruction images.

Table 1 shows the reconstruction errors. From this table, we can see that the original images can be perfectly reconstructed by MagnitudeCut method when *s_n* is equal to the triple of the size of original images.

Table 1 The reconstructed error of the MagnitudeCut algorithm according to the Fig. 2.

| **Image** | **(a)** | **(b)** | **(c)** | **(d)** | **(e)** | **(f)** |
|---|---|---|---|---|---|---|
| **Error** | 2.0145e-16 | 4.1772e-14 | 0.2212e-12 | 8.1321e-14 | 3.8129e-15 | 9.5183e-15 |

### *3.3. 2-D complex symmetric images reconstruction*

In some cases, we have to deal with the symmetric image. Whereas, the methods that mentioned in the introduction under conditions proposed by Oppenheim and Lim [13], Levi and Stark [2] cannot recover the symmetric image. The original image $\boldsymbol{x} \in \mathbb{C} \times \mathbb{C}$, and suppose we know the phase of its Fourier transform. The reconstructed image is shown in Fig. 3. Table 2 shows the reconstruction errors.As before, we can see that the symmetric image can also be recovered from the phase of its Fourier transform by the proposed MagnitudeCut method.

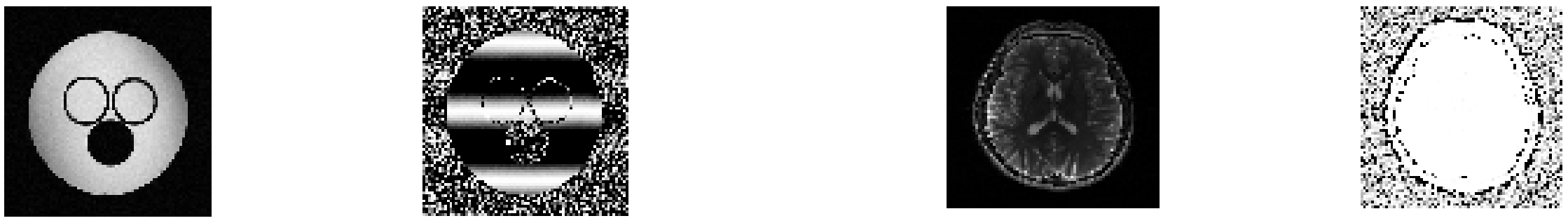

(a) (b) (c) (d)

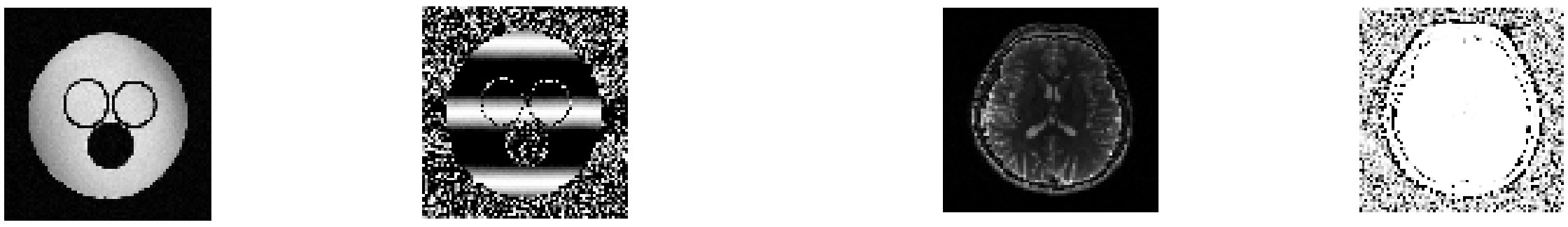

(A) (B) (C) (D)

Fig. 3 Water Phantom (left) and Brain Image (right). (a) and (c) are the magnitude of the original images, (b) and (d) are the phase of the original images. (A) and (C) are the magnitude of the reconstruction images, (B) and (D) are the phase of the reconstruction images.

Table 2 The reconstructed error of the MagnitudeCut algorithm according to the Fig 3.

| **Image** | **Water Phantom** | **Brain Image** |
|---|---|---|
| **Error** | 5.0465e-15 | 6.1168e-14 |

### *3.4. Comparison with other methods*

In this section, we compare the reconstructed results using the proposed method with the Greedy algorithm [18], and Iterative algorithm [13]. Assume the original signal is $\boldsymbol{x} \in \mathbb{R}^{32}$, and the phase of its Fourier transform is given. The range of iterations of these methods was set between 32 and 320. Fig. 4 exhibits the results. In this figure, we list the cost time, the reconstruction errors ($\left\|\boldsymbol{x}-\hat{\boldsymbol{x}}\right\|_2^2 / \left\|\boldsymbol{x}\right\|_2^2$) and the variance of the error under different sampling numbers by the three methods, where the symbol *s_n* denotes the sampling number. To the second and third column images of Fig. 4, for convenient observation, we set the vertical coordinates to the logarithmic function ( log10() ) value of the reconstruction error and the variance of the error. Note that for accurate comparison between MagnitudeCut and Iterative algorithm, Greedy algorithm, we take an average of errors and variances of experiments repeated 1000 times.

As we can see, with the increase of *s_n*, all the three methods can completely reconstruct the original signal. However, compared to the Iterative algorithm and Greedy algorithm that need at least three times that of the number of the original signal, our method needs only two times. What's more, when *s_n* is equal to the number of the original, the MagnitudeCut method can still achieve better recovery performance than the other two methods. But beyond that, we find that as the number of phase information changes, the reconstruction error when using MagnitudeCut varies at a smaller range than the other two methods from the third column images of the Fig. 4, which means the MagnitudeCut is more stable. Although the first image of the third column of Fig. 4 shows the Greedy method is stable, it cannot reconstruct the original signal.

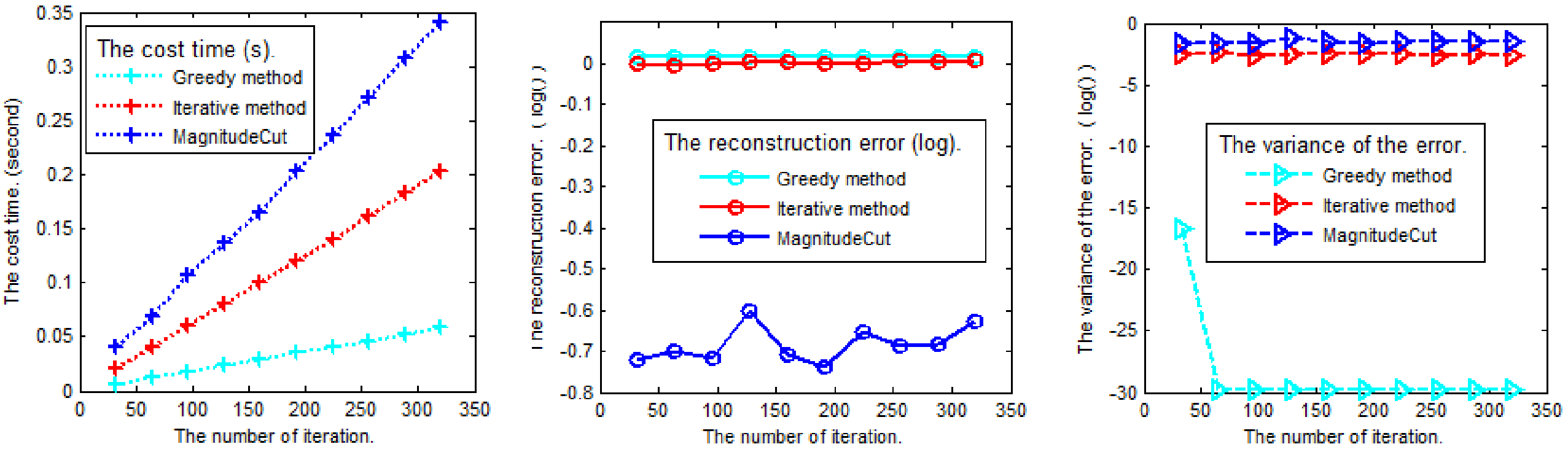

(a) The sampling number *s_n* is equal to 32.

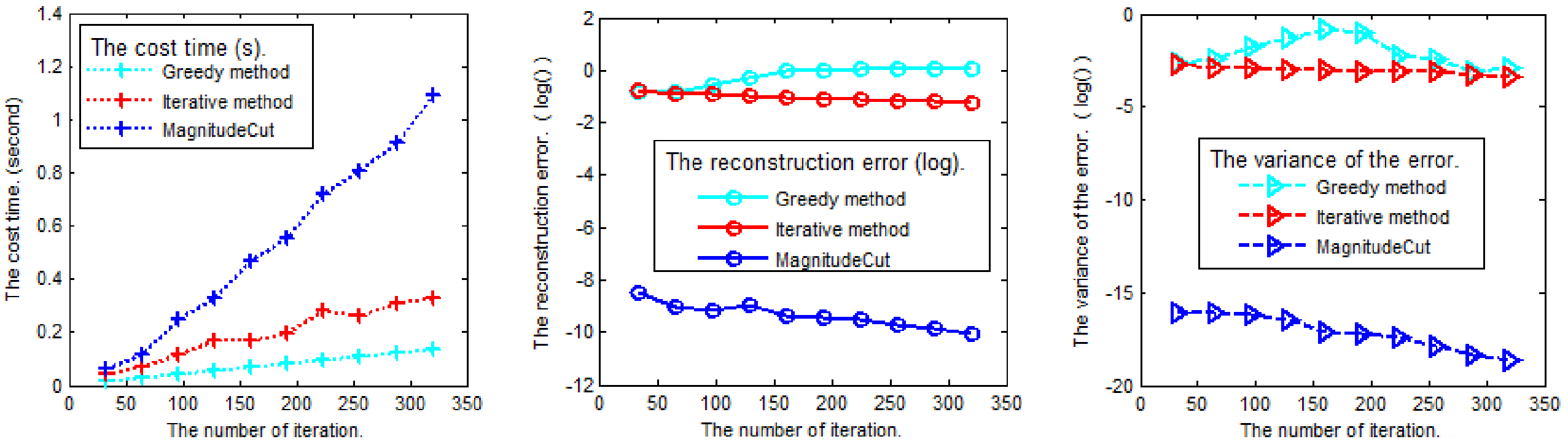

(b) The sampling number *s_n* is equal to 64.

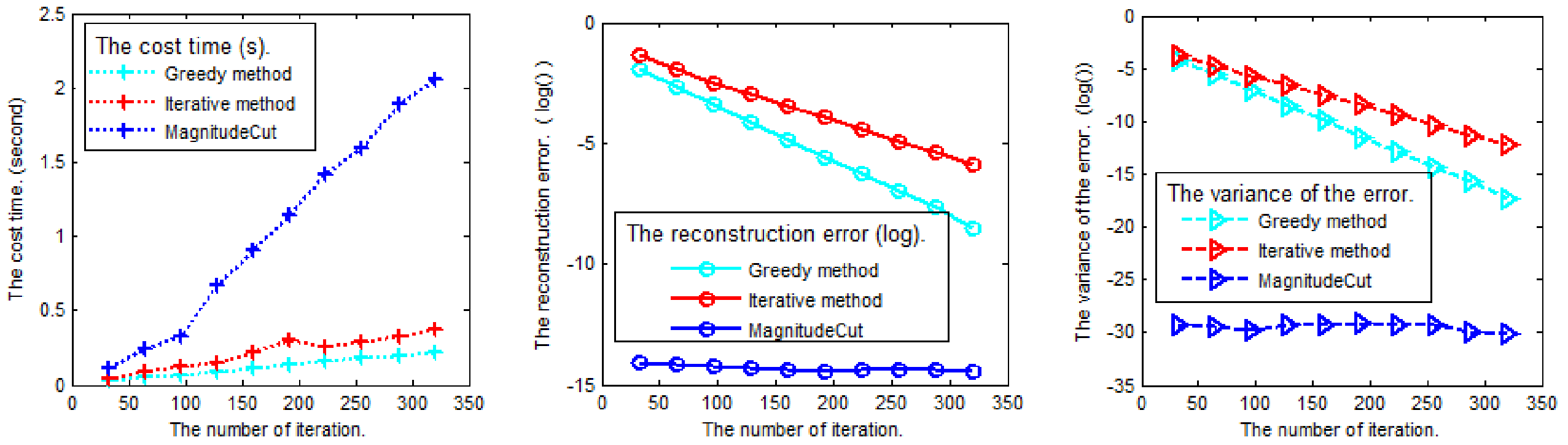

(c) The sampling number *s_n* is equal to 96.

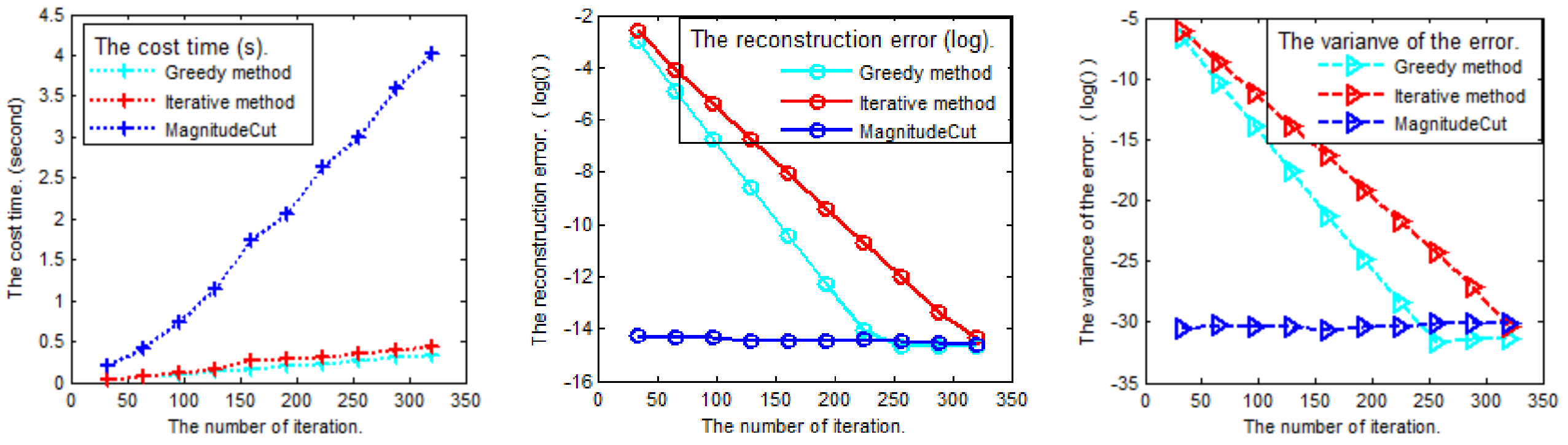

(d) The sampling number *s_n* is equal to 128.

Fig.4 The cost time, reconstructed error and the variance of the error of the MagnitudeCut algorithm vs the traditional approaches under different sampling numbers.

The reason for the stability of the Greedy method is that it may quickly get a local optimal solution on the condition that the sample number is little. From the first column images of Fig. 4, we can see that the proposed method requires more computation time compared to the other two methods. However, it is important to note that with about 32 iterations, our method can already achieve significantly better results than the other two methods with basically the same computation time.

## 4. Conclusion and future work

In this paper, we have proposed a new algorithm for solving the problem of signal reconstruction from the phase only information. Numerical results have been provided to illustrate the feasibility of the algorithm. The merit of the proposed algorithm is that it can also reconstruct the symmetric image from its Fourier phase and less sampling number of the phase information is needed to reconstruct the original signal. Phase information can preserve many important features of a signal and the magnitude information is distorted in some cases, for example, in long-term exposure to atmosphere turbulence or when images are blurred by severely defocused lenses with circular aperture stops. Therefore, if the phase information is used to describe the signal features, the requirements for the storage and the transmission bandwidth can be reduced. Just as the PhaseCut algorithm [12] is the basis of the scattering convolution networks [21], our algorithm shows that we may also construct a new scattering convolution network by using only the phase information. In the next work, we will consider quaternion signal. What's more, we will extend the MagnitudeCut algorithm for practical signal with noise.

## Acknowledgments

This work was supported by the National Natural Science Foundation of China under Grants 61201344, 61271312,61401085, 81101104, 11301074 and 61073138, by the Research Fund for the Doctoral Program of Higher Education (No. 20120092120036), the Project-sponsored by SRF for ROCS, SEM, and by Natural Science Foundation of Jiangsu Province under Grants BK20150650, DZXX-031, BY2014127-11, by the '333' project under Grant BRA2015288, by the Qing Lan Project, and by the High-end Foreign Experts Recruitment Program (GDT20153200043). The authors are thankful to Dr. Xiaobo Qu of Xiamen University, China, for providing the Water Phantom image and Brain Image in Fig. 3.